%% file: main.tex
\documentclass[pagesize]{scrartcl}
\usepackage{mystyle}
\pgfplotsset{compat=1.18}

\definecolor{frenchrose}{rgb}{0.96, 0.29, 0.54}

\newcommand{\inner}[2]{\left\langle #1, #2 \right\rangle}
\title{A First Step Towards Mesh-Free \\ Probabilistic Shape Optimization}
\author{
    Stephan Schmidt\thanks{Trier University, Department IV -- Mathematics. M.\ W\"{u}rschmidt gratefully acknowledges financial support from the German Research Foundation (DFG) within the Research Training Group 2126: Algorithmic Optimization.} 
    \and 
    Maximilian W\"{u}rschmidt\footnotemark[1]
}
\date{}
\begin{document}

\allowdisplaybreaks
\maketitle

\begin{abstract}
	We present an initial implementation of a probabilistic PDE-constrained shape optimization algorithm. Our method is based on a novel probabilistic representation of the shape derivative, which is evaluated using Monte Carlo sampling; and does not rely on a mesh. The underlying state is represented with a neural network-based PDE solver on point clouds.
	The methodology is applied throughout to a benchmark tracking problem.

    %
    
\end{abstract}

\section{Introduction}
%
Tracking type shape optimization governed by elliptic partial differential equations are prototypi\-cal for a wide class of mathematical inverse problems. They are of relevance for numerous applications, in particular electrical and geological impedance tomography. In this class of problems, one seeks to reconstruct an inclusion indirectly by observing potentials. The inclusion effects the conductivity of the base material and the optimizer has to find material conductivity parameters, such that a simulated potential (state) coincides with the measured potential (tracking data), see, e.g., \cite{doi:10.1137/S0036144598333613}. A sharp interface between inclusion and base material can be realized via two subdomains. The problem is to optimize the shape of the interface between the two. If the inclusion is a proper void, it is often sufficient to consider just a single domain. 
In its most general form, see, e.g., \cite[Section~2.6]{Sokolowski_Shape}, this leads to the following problem
%
\begin{equation}\label{eqn:TrackingProblem}
    \min_{\Omega \subset \HoldAll } \int_{\Omega} \tfrac{1}{2} \big( u(x) - \trackingData (x) \big)^2 \ \lambda^2 (\de x) 
    \quad \text{s.t.}\begin{cases}
                        -\Delta u &= 1 \quad \text{on }\Omega ,\\
                        \hspace*{0.55cm} u &= 0 \quad\text{on }\partial\Omega. 
                    \end{cases} 
\end{equation}
where $\HoldAll \subset \R^2$ is a rectangular \textit{hold-all} domain and $z_d\colon\HoldAll \to \R$ is the \textit{tracking data} map.\footnote{ 
	In this context, the shape functional is
	$J \colon (u, \Omega) \to \tfrac{1}{2} \|u- \trackingData \|_{L^2(\Omega)}^2$
	where $(u, \Omega)$ is to be understood as a one-to-one correspondence, i.e., we assume that for each admissible shape $\Omega$ the PDE-constraint has a unique classical solution $u$. E.g., this is the case whenever an exterior sphere condition is satisfied, see \cite[Theorem~6.13]{GilbargTrudinger}.
}

When the data is noisy such inverse problems usually require an additional regularization term. Typically this functional is the volume or perimeter of the domain, fostering the reconstruction of smooth domains. Nevertheless, some applications, e.g., geological impedance problems, demand reconstructing a non-smooth domain, see, e.g., \cite{TVGeoInverse_discrete, TVGeoInverse_continuous}.

In gradient based shape optimization an (locally) optimal shape can be approximated with Algorithm~\ref{algorithm:ShapeGD}. Subsequently, we present a novel methodology in which we use a neural PDE solver to compute the state and we employ a probabilistic representation of the shape derivative to compute the shape gradient.
\begin{algorithm}
	\SetAlgoLined
	\textbf{initialization} Choose initial domain $\Omega$ and tolerance $\varepsilon_{\text{TOL}}$\\
	\While{$\| \nabla J (\Omega) \| > \varepsilon_{\text{TOL}}$}{
		1. Calculate state $u$ on $\Omega$ \\
		2. Calculate shape gradient $\nabla J (\Omega)$ \\ 
		3. Perform linesearch to obtain $\alpha$ [optional]\\
		4. Apply geometric flow of $\alpha \nabla J (\Omega)$ to $\Omega$
	}
	\caption{Shape Gradient Descent}\label{algorithm:ShapeGD}
\end{algorithm}

In recent years, deep learning based algorithms, such as so called physics-informed neural networks (PINNs), as well as, exploiting the stochastic setting, different variations the deep BSDE solver have been proposed, analyzed and are indeed capable of approximating PDE solutions. For an overview we mention \cite{beck2023overview, NuskenRichter2023, raissi2019physics, tanyu2023deep} and \cite[Part IV]{jentzen2023mathematical}. Furthermore, these algorithms have the potential to overcome the curse of dimensionality, see, e.g., \cite{GrohsHerrmann2022}.

The probabilistic representation of the shape derivative was developed in \cite{schlegel2024probabilistic} for a general state space dimension. In this paper we consider a $2$D benchmark problem, nevertheless we emphasize that the probabilistic formulation as well as the below discussed (meshless) Monte Carlo simulations do not suffer from the curse of dimensionality, see, e.g.,
\cite{caflisch1998monte, gelman1997weak}.

Thus, a key advantage of this approach is that the current shape is represented solely by a point cloud or an equivalent level-set description, making the framework amenable to extension to high-dimensional settings.

%
At this stage, our implementation requires a boundary mesh, primarily for the computation of the shape gradient (Step 2 in Algorithm~\ref{algorithm:ShapeGD}) and, for convenience, in parts of the evaluation of the probabilistic shape derivative.
\footnote{The code is publicly available at {\scriptsize \url{https://github.com/max-wuer/ProbabilisticShape}}}
%

\section{Calculate State via PINN}\label{sec:statePINNs}

We approximate the state using a physics-informed neural network (PINN), see \cite{raissi2019physics, Sirignano2018}. This is an unsupervised, resp. self-contained, learning algorithm -- we refer to \cite[Chapter 16]{jentzen2023mathematical} for a comprehensive introduction -- applied to feedforward neural networks $\nn \defined \nn^\theta \colon \R^d \to \R$ with smooth activation function. Here $\theta$ denotes the networks weights and is given by a tuple of weight matrices and bias vectors -- we refer to \cite[Chapter 2]{petersen2024mathematical} for an in-depth introduction. We remark that once an activation function is specified, each network is uniquely determined by its weights.    
%

On a given bounded domain $\Omega \subset \HoldAll$ we draw collocation points $\Xdomain\times\Xbdry \subset \Omega\times\partial\Omega$ according to a predetermined distribution. The objective in the training process is to compute a minimizer of the PINN loss functional, which, for the Poisson equation corresponding to problem \eqref{eqn:TrackingProblem}, is given by 
\begin{equation}\label{eqn:PINNLoss}
	\lossPINN\colon \nn \mapsto \tfrac{1}{|\Xbdry|+|\Xdomain|} \Big( 2|\Xbdry| \sum_{x \in \Xbdry } \big| \nn(x) \big|^2 
	+ |\Xdomain| \sum_{x \in \Xdomain } \big| \Delta \nn(x) + 1 \big|^2 \Big)
\end{equation}
Note the PINN functional is depending on the sample of collocation points $\Xdomain\times\Xbdry$. Thus
\begin{equation}
	\lossPINN (\nn) = \lossPINN (\nn;\ \Xdomain\times\Xbdry) 
	= \lossPINN (\nn^\theta;\ \Xdomain\times\Xbdry) 
	= \lossPINN (\theta;\ \Xdomain\times\Xbdry) 
\end{equation}
Assuming that the collocation points are drawn independently from a continuous distribution with full support on $\Omega$, the law of large numbers implies that any dependence induced by the sampling procedure becomes asymptotically negligible.

The minimization of $\lossPINN$ is accomplished by gradient descent, see, e.g., {\cite[Chapter~10]{petersen2024mathematical}}. 
Here the computation of the Laplacian $\Delta \nn$ and of PINN loss gradient $\Diff_\theta \lossPINN$ are accomplished via automatic differentiation.

\section{Calculate Shape Gradient via Probabilistic Shape Derivative}\label{sec:DeformationEQ}
%
In this section, we discuss the second step of Algorithm~\ref{algorithm:ShapeGD}, which consists of computing the shape update, i.e., the {effective shape deformation for each descent step}. As in Section~\ref{sec:statePINNs}, our presentation considers a generic \textit{current domain} $\Omega\subset\HoldAll$ and explains how to compute the shape update.

This discussion is divided into two parts: first, we present the probabilistic representation of the shape derivative from \cite{schlegel2024probabilistic}, which is based on a Feynman-Kac expansion of the boundary sensitivity of the underlying state and expresses the derivative in terms of expectations of random variables that can be simulated without a mesh; second, we address the computation of the shape gradient, which at this stage requires us to introduce a boundary mesh.

\begin{remark}\label{remark:NotationShapeDerivative}
	For readers with a background in shape optimization, we emphasize that the probabilistic representation $\D J$, see \eqref{eqn:ProbShapeFuncDerivative} below, is based on a pre-image/pull-back perturbation of identity. Hence the probabilistic approach yields the classically considered directional shape derivative $\mathcal{D} J$, based on push-forward perturbation of identity, with flipped sign, i.e., 
	\[
	\D J [\dir] = \mathcal{D} J[-\dir] \qquad\text{for all admissible perturbations }\dir\closeEqn
	\]
\end{remark}

For a comprehensive introduction to shape optimization we refer to the recent exposition \cite{Alllaire2021}. Particularly for the classical notion of the push-forward shape derivative, see \cite[Definition~4.1]{Alllaire2021}.

\begin{remark}\label{remark:NotationDomainBoundaryContVSDisc}
	We use the notation $\partial\Omega$ for the meshless boundary and $\partial\Omega_{\text{FEM}}$ for the respective boundary mesh.
	Moreover we write $\Omega_{\nn}$ whenever a PINN approximation of the state is used to describe the current domain numerically. \close
\end{remark}

\paragraph*{Probabilistic Shape Derivative} 
For an admissible perturbation $\dir\colon\partial\Omega\to\R^2$ the formula\-tion in \cite[Theorem 3.7]{schlegel2024probabilistic} applied to our tracking type problem reads as follows
\begin{multline}\label{eqn:ProbShapeFuncDerivative}
	\D J [\dir] = \constant{+} \E\Big[\big\langle \dir,\nabla u \big\rangle \big(\enX^{X_0^+}\big)\Big] -\constant{-}\E\Big[\big\langle\dir ,\nabla u \big\rangle\big(\enX^{X_0^-}\big)\Big]\\
	- \int_{\partial\Omega} \tfrac{1}{2} \trackingData(y)^2\,\big\langle\dir,n\big\rangle(y)\,\dSphere{1}{y}
\end{multline}
where $\enX^{X_0^+}$, $\enX^{X_0^-}$ are so-called random-start exit-kill random variables taking values in $\partial\Omega$ and $m^{+}, m^{-}\geq 0$ are constants, which will be introduced below. 
	
Both expectations in \eqref{eqn:ProbShapeFuncDerivative} have similar structure, merely the support of the random-start distribution differs. 
For the tracking type problem, given the current state $u$, consider the partition\footnote{
	Notice in terms of the objective in \eqref{eqn:TrackingProblem} the partition corresponds to level set $\Omega^+ = \{ x\in\Omega\ |\ \partial_y \phi (x, u(x)) \geq 0 \}$ of the map $\phi\colon (x,y)\mapsto \tfrac{1}{2} ( y - \trackingData (x) \big)^2$.
} of the domain 
\[
	\Omega^+ \defined \bigl\{x\in\Omega \,\big|\, u(x)\geq \trackingData(x)\bigr\}
	\qquad\text{and}\qquad
	\Omega^- \defined \Omega\setminus\Omega^+
\]
We use $\pm$ as a placeholder to indicate dependence on the chosen subdomain. Accordingly
\begin{equation}\label{eqn:ProbShapeConstants}
	\constant{\pm} = \pm\int_{\Omega^\pm} u(x)-\trackingData(x)\dLebesgue{2}{x}
\end{equation}
The random-start exit-kill variables consist of a random start $X_0^\pm$, which satisfy the distribution $\rho^\pm \dLeb{2}$ with explicitly available densities
\begin{equation}\label{eqn:RandomStartDensities}
	\rho^\pm\colon\R^d \to [0, \infty); \
	x \mapsto 
	\tfrac{1}{\pm\constant{\pm}} \big(u(x)-\trackingData(x)\big)\
	\1_{ \Omega^\pm } (x)
\end{equation}
that depend on the partition of $\Omega$. 
For the tracking problem killing is absent and thus, the exit is characterized by the dynamics of a Brownian motion $W$, i.e.,
the stopping time $\tau \defined \inf \{ t\geq 0 \ |\  X_0^\pm + \sqrt{2} W_{t} \notin \Omega \}$. The random variables in \eqref{eqn:ProbShapeFuncDerivative} are then given by the above described process\footnote{
	Note this is a continuous-time process, stopped at first exit from an open set, and hence, attains values in $\partial\Omega$.
}
at first exit;
\begin{equation}\label{eqn:RandomStartExitKill}
	\enX^{X_0^\pm} = X_0^\pm + \sqrt{2} W_{\tau}
\end{equation}

%
%
The remainder of this paragraph is devoted to a presentation of a mesh-free evaluation of the probabilistic representation of the shape derivative, see \eqref{eqn:ProbShapeFuncDerivative}. 

Throughout, we assume that $\nn \approx u$ is a trained PINN for the state on the current domain $\Omega$, i.e., $\nn$ is a good approximation of a (local) minimizer of the PINN loss \eqref{eqn:PINNLoss}. 
We employ the PINN to approximate $\Omega \approx \Omega_{\nn} \defined  \{\nn >0\} $, which is justified by the maximum-principle respectively the Hopf-Lemma, see \cite{hopf1952remark}. Similarly, we approximate the domain partition with 
\[
	\Omega^+_{\nn} \defined \{\nn\geq\trackingData\}\cap\{\nn\geq 0\} \qquad\text{and}\qquad \Omega^-_{\nn} \defined \{\nn\leq\trackingData\}\cap\{\nn\geq 0\}
\]

The simulation of the random-start and the evaluation of the constants is obtained using  acceptance-rejection sampling, see Algorithm~\ref{algorithm:AccRej}. 
We refer to, e.g., \cite[Section~2.2.2]{Glasserman} for further details. 
In all applications we use a uniform distribution on the hold-all domain as respective reference density, i.e., $\rhoREF = \tfrac{1}{\lambda^2 (\HoldAll)} \1_{\HoldAll}$.

%
%
\begin{algorithm}
	\SetAlgoLined
	\textbf{requirement} Reference density $\rhoREF\colon\R^d\to [0,\infty)$ from which we can sample and such that there is $C>0$ satisfying $\rho (x) \leq C \rhoREF(x)$ for all $x\in\R^d$.\\
	\textbf{initialization} Choose number of samples $M\in\nat$ to be drawn {independently} from target density $\rho$. Set $\textbf{X}^\rho = \emptyset$. \\
	\While{$ | \textbf{X}^\rho  | < M $}{
		1. Sample $X\sim \prob_{\text{ref}} = \rhoREF \dLeb{d}$ \\
		2. Sample $U\sim \text{Uniform} (0,1) $ \\ 
		\If{ $U \leq \tfrac{\rho(X)}{C\rhoREF(X)}$ }{ Accept sample: $\textbf{X}^\rho \gets \textbf{X}^\rho \cup \{X\} $ }
	}
	\caption{(Batch) Acceptance Rejection Sampling}\label{algorithm:AccRej}
	\textbf{return} Set of samples $\textbf{X}^\rho$. 
\end{algorithm}
%
%
Note for a uniformly distributed $\mathcal{U}^\pm$ on $\Omega^\pm$ it holds that
\begin{equation}\label{eqn:ProbShapeConstantDecomposition}
		\constant{\pm} 
		\approx \pm \int_{\Omega^\pm} \nn (x)-\trackingData(x) \dLebesgue{2}{x}
		= \pm \lambda^2 (\Omega^\pm)\, \E\big[ \nn (\mathcal{U}^\pm)-\trackingData(\mathcal{U}^\pm) \big]
\end{equation}
We approximate each factor individually using Algorithm~\ref{algorithm:AccRej}. For both we use the reference density $\rhoREF$ from above, the constant $C=1$ and as target the density $\rho=\tfrac{1}{\lambda^2 (\Omega^\pm_{\nn})} \1_{\Omega^\pm_{\nn}}$, to generate a set of samples $ \samples{\Omega^\pm_{\nn} }{uniform} $. 

Note $\lambda^2 (\Omega^\pm) = \lambda^2 (\HoldAll)\ \prob_{\text{ref}} [ \Omega^\pm ]$. 
Thus, the Lebesgue measure is obtained as a special case, where the Monte Carlo estimate reduces to the ratio of accepted samples to the total number of trials, i.e.,
\begin{equation}\label{eqn:ApproxMeasurePartition}
	\prob_{\text{ref}} [{\Omega^\pm}] \approx 
	\prob_{\text{ref}} [{\Omega^\pm_{\nn}}] \approx
	\frac{ \# \text{ successful while loops} }{  \# \text{ total while loops} }
\end{equation}
Particularly, the specific coordinates of the samples $\samples{\Omega^\pm_{\nn} }{uniform}$ are not used here. Next, with the same samples, we evaluate the second factor, i.e., the expectation
\[
	\E\big[ \nn (\mathcal{U}^\pm)-\trackingData(\mathcal{U}^\pm) \big]
	\approx \tfrac{1}{|\samples{\Omega^\pm_{\nn} }{uniform}| }  \sum_{x\in\samples{\Omega^\pm_{\nn} }{uniform}} \nn(x) - z_d (x) .
\]
With those approximations at hand, we set $\constantNN{\pm}$ to denote the respective numerically computed versions of \eqref{eqn:ProbShapeConstantDecomposition}.

The random-start exit-kill variables \eqref{eqn:RandomStartExitKill} are sampled sequentially; specifically, Algorithm~\ref{algorithm:AccRej} generates sets of random starts $\samples{0}{$\pm$}$, from which the exit coordinates on the domain boundary are computed using an Euler-Maruyama scheme, see Algorithm~\ref{algorithm:RandomStartExitKill}.

To this end, we invoke the acceptance-rejection sampling, see Algorithm~\ref{algorithm:AccRej}, with the uniform distribution on the hold-all domain as reference density $\rhoREF$. The target densities are given as explicit mappings $\rho^+, \rho^-$ in \eqref{eqn:RandomStartDensities} and approximated using the PINN level set, i.e.,
\[
	\rho^\pm_{\nn} \defined  
	\tfrac{1}{\constantNN{\pm}} \big| \nn -\trackingData \big|\ \1_{ \Omega^\pm_{\nn} } 
\]
and as constants for the acceptance condition (we use the same samples $\samples{\Omega^\pm_{\nn}}{uniform}$ as above);
\begin{equation}\label{eqn:acc-ref-constantRandomStart}
	C^\pm = \lambda^2(\HoldAll)\ \safety \ \max \big\{ \rho^\pm_{\nn} (x)\ \big|\ x\in \samples{\Omega_{\nn}^\pm}{uniform} \big\}
\end{equation}
where $\safety >1$ denotes a safety factor, to ensure that the estimate in the algorithm can be reasonably extrapolated from the uniformly drawn point clouds to the entire level set $\Omega_{\nn}^\pm$. From this we obtain the random-starts $\samples{0}{$\pm$}$.
In a second step, we compute the corresponding exit points $\exitSamples $, as presented in Algorithm~\ref{algorithm:RandomStartExitKill}. Here only the projection onto $\partial\Omega_{\text{FEM}}$ requires the underlying boundary mesh. The reason is twofold: on the one hand, it is a commodity to use the mesh-normals and on the other hand, we require to hit the boundary with machine precision since finite elements are not robust for evaluations outside their underlying mesh.\footnote{
	For completeness we mention that at each vertex of $\partial\Omega_{\text{FEM}}$ is a cone without well-defined normal to project along. In case the first exit lies in said cone, we project onto the respective vertex. Overall we always project onto the closest point in $\partial\Omega_{\text{FEM}}$.
	Specifically, if $x\in\R\setminus\Omega_{\nn}$ is a realization of a first Euler-Maruyama exit, let $v$ denote the closest and $w$ denote the second-closest vertex to $x$ in $\partial\Omega_{\text{FEM}}$. Correspondingly, let $n_{\text{FEM}}$ be the outer normal of the edge $(v,w)$, we then set
	$ x_p \gets (\mathcal{I} - n_{\text{FEM}}n_{\text{FEM}}^\top) x + n_{\text{FEM}}n_{\text{FEM}}^\top v $ to obtain the closest point on $\partial\Omega_{\text{FEM}}$.
}
%
%
\begin{algorithm}
	\SetAlgoLined
	\textbf{initialization} Sample start $x_0 \gets X_0^\pm$. Set discrete-time stepsize $\Delta \in (0,1)$.\\
	\While{$ x_k \in \Omega_{\nn} $}{
		Euler Maruyama Update: $x_k \gets x_{k-1} + \sqrt{2} \sqrt{\Delta} \xi$ for $\xi\sim\mathcal{N} (0,\mathcal{I})$
	}
	\textbf{return} 
	Orthogonal projection of $x_k$ onto boundary $\partial\Omega_{\text{FEM}}$
	\caption{Random-Start Exit(-Kill) Variables Sampling}\label{algorithm:RandomStartExitKill}
\end{algorithm}
%
%

The approximation of the boundary integral in \eqref{eqn:ProbShapeFuncDerivative} uses the boundary mesh $\partial\Omega_{\text{FEM}}$ and the corresponding facet normals. At this point, 
$\partial\Omega_{\text{FEM}}$ is employed out of convenience since it is already available.
Alternatively, one could employ a uniform distribution on $\partial\Omega$ and perform a Monte Carlo simulation as above. In that case the reference measure is the $1$D Hausdorff measure restricted to $\partial\Omega$, i.e., $\mathcal{S}^1 (\cdot\cap\partial\Omega)$. The respective approximations read as follows
\begin{equation}\label{eqn:ProbShapeStokes}
	\int_{\partial\Omega_{\text{FEM}}} \!\! \tfrac{1}{2} \! \big[ \trackingData^2 \langle\dir, n_{\text{FEM}}\rangle \big] (y) \dSphere{1}{y}
	\approx \int_{\partial\Omega} \!\! \tfrac{1}{2} \! \big[ \trackingData^2 \langle\dir, n\rangle \big] (y) \dSphere{1}{y}
	= \prob_{\mathcal{S}^1} (\partial\Omega)\ \E \Big[ \tfrac{1}{2} \big[ \trackingData^2 \langle\dir, n\rangle \big] (\mathcal{U}) \Big]
\end{equation}
where $\mathcal{U}$ is uniformly distributed with respect to the measure 
$\prob_{\mathcal{S}^1} \defined \tfrac{1}{\mathcal{S}^1 (\partial\Omega)} \mathcal{S}^1 (\cdot\cap\partial\Omega)$. 

%
%


\paragraph*{Deformation Equation}
The probabilistic shape derivative~\eqref{eqn:ProbShapeFuncDerivative} is providing us with a directional derivative, which we turn into a gradient via the Riesz map
\begin{equation*}
 \inner{W}{\dir}_{H^1(\partial \Omega)} \stackrel{!}{=} \D J [\dir] \quad \text{for all admissible directions } \dir.
\end{equation*}
Using the $H^1$-inner product has several advantages. On the one hand, it ensures that the unknown domain boundary $\partial \Omega$ stays within the desired regularity class. On the other hand, this procedure can be interpreted as an approximate Newton scheme, where the $H^1$-inner product is an approximation to the Hessian of the problem. For more details see~\cite[Chapter~3]{Schmidt-Schulz-2023}, where the actual Hessian is used. While there are mesh-free implementations of the Laplace-Beltrami-Operator, i.e., the PDE of the $H^1$-inner product, see, e.g., \cite{SUCHDE20192789, TAMAI201779}, these methods are typically non-trivial and outside the scope of this work.

We now introduce a boundary mesh $\partial\Omega_{\text{FEM}} = (\{x_j\}_{j\leq k}, \mathcal{E})$, i.e., a closed polygonal chain in 2D with vertices $\{x_j\}_{j\leq k}$, edges $\mathcal{E}$ and an outer normal $n_\text{FEM}$, constant on each edge, see also Remark~\ref{remark:NotationDomainBoundaryContVSDisc}. On this mesh, we declare the first order continuous Lagrange finite element space;
\begin{equation*}
	\mathcal{V}_\text{FEM} \defined \big\{v \in \C(\partial \Omega_\text{FEM}, \R) \ \big| \ {v} |_{e} \in P_1(e) \ \forall \ e \in \mathcal{E} \big\}.
\end{equation*}
Here, $P_1$ is the set of all first order polynomials and $\mathcal{E}$ is the set of edges of the mesh.
We equip the space $\mathcal{V}_\text{FEM}$ with the base $\{v_j\}_{j\leq k}$ with the property
\begin{equation*}
	v_j(x_i) = \delta_{ij} \defined \begin{cases}1, &\quad i=j\\0, &\quad i\neq j\end{cases}
\end{equation*}
Each $v_j$ has its support on exactly two neighboring edges $e^\ell$, $e^r$ with respective constant outer normals $n_\text{FEM}^\ell$, $n_\text{FEM}^r$.
As such, we can define our set of admissible deformations as the span of
\begin{equation*}
	V_j \colon \partial\Omega_{\text{FEM}}\setminus \left(\partial e^\ell \cup \partial e^r\right) \to \R^2;\quad  x\mapsto  \begin{cases} v_j(x)n_\text{FEM}^\ell, &\quad x \in e^{\ell, \circ}\\ v_j(x)n_\text{FEM}^r, &\quad x \in e^{r,\circ}\\ 0, &\quad\text{otherwise}\end{cases}
\end{equation*}
where $e^\circ$ denotes the interior of the respective edge. A definition on the (three) vertices $\partial e^\ell \cup \partial e^r$ is not necessary because these points never happen to be exit points in the probabilistic shape derivative. Indeed, the $H^1(\partial \Omega_\text{FEM})$-gradient of the probabilistic shape derivative~\eqref{eqn:ProbShapeFuncDerivative} can now be determined using classical finite elements, i.e., by finding $w \in \mathcal{V}_\text{FEM}$ satisfying
\begin{equation*}
	\sum_{e\in \mathcal{E}} \sum_{s=1}^{k_e} \omega_j^e\left(\frac{1}{2}\inner{\nabla_\Gamma {{v_j}|_{e}(\tilde{x}_s^e)}}{\nabla_\Gamma w|_{e}(\tilde{x}_s^e)} + \inner{{v_j}|_{e}(\tilde{x}_s^e)}{w|_{e}(\tilde{x}_s^e)} \right)
	\stackrel{!}{=} \D J [V_i] \quad \forall i = 1,\dots, k.
\end{equation*}
Here, $\nabla_\Gamma$ denotes the tangential gradient\footnote{Note this is the natural gradient concept for function in $\mathcal{V}_\text{FEM}$} and $\omega_j^e$ is the quadrature weight at $\tilde{x}_s^e$ on the edge $e$ with $k_e$ quadrature positions. To update the shape using this $H^1(\partial\Omega_\text{FEM})$-gradient, we then compute the Galerkin projection $W$ of $w\,n_\text{FEM}$ into the space $\mathcal{V}_\text{FEM} \times \mathcal{V}_\text{FEM}$ and set the next domain to be
\begin{equation*}
	\partial \Omega_{\text{FEM}}^{\ell + 1} = \left\{ x_j + \tau W(x_j) \ \vert \ x_j \text{ vertex in } \partial \Omega_\text{FEM}^{\ell}\right\},
\end{equation*}
where $\alpha > 0$ is a suitable step length. Recall Remark~\ref{remark:NotationShapeDerivative} in regard to flipped the sign.

\section{Numerical Results}\label{sec:numerics}

In the simulation (see Figure~\ref{fig:SimulationResult_ATL_ENUMATH}) we consider a case where the optimal shape is explicitly available. Namely we set $\trackingData$ to be the solution of the PDE-constraint with the domain $\Omega^*$ being an ellipsis centered at $c=(\tfrac{1}{2},\tfrac{1}{2})$, with semi-major axis $a=0.4$ and semi-minor axis $b=0.3$; extended with constant zero into the hold-all $\HoldAll$.\footnote{
	Note this map is explicitly available: $\trackingData\colon \HoldAll \to \R;\ (x_1,x_2) \mapsto \tfrac{a^2b^2}{2(a^2+b^2)} \big( 1- \tfrac{(x_1-c_1)^2}{a^2} -\tfrac{(x_2-c_2)^2}{b^2} \big) \ \1_{\Omega^*}(x_1,x_2)$.
}
In that case the ellipsis $\Omega^*$ solves \eqref{eqn:TrackingProblem}.

\input{optimization_results_figure}

We initiate the shape gradient-descent scheme, see Algorithm~\ref{algorithm:ShapeGD}, with a ball $B_R (c)$ with the same center as the ellipsis and radius $R=0.25$. This framework allows to set the hold-all domain to $\HoldAll = [0,1]^2$, and thus the factor $\lambda^2 (\HoldAll) = 1$ is negligible.

%
%
We employ a fully-connected feedforward neural network with two hidden layers of $20$ neurons each and $\tanh$ as activation function. Particularly, we approximate the PDE on the hold-all
domain $\HoldAll$ with an interior Dirichlet interface at $\partial\Omega$. 
The main reason is that we use the PINN to describe the current domain via the level-set $\{\nn\geq 0\}$. Moreover, this approach considerably improves the approximation quality of the state gradient.

The underlying collocation points $\Xdomain\times\Xbdry$ consist of point clouds sampled independently. The interior points $\Xdomain$ are drawn according to a uniform distribution on $\HoldAll$. The Dirichlet points $\Xbdry$ are generated using the boundary mesh $\partial\Omega_{\text{FEM}}$ by sampling uniformly from each edge individually.
Specifically we use {$|\Xdomain| = 45000$ and $|\Xbdry| = 22500$} samples. For the Dirichlet points, the same number of samples is drawn from each edge, regardless of its geometric length.

The PINN optimization is performed using the \textsc{Adam} optimizer, see \cite{kingma2015adam}, with a piecewise constant learning rate. The implementation is carried out using the \textsc{PyTorch} library, see \cite{paszke2019pytorch}. 
We emphasize the large difference in the number of optimization steps: the initial training is long, requiring {$50000$ -- $100000$} gradient descent steps and has been restarted multiple times to find a very good local optimum.
In each subsequent iteration of Algorithm~\ref{algorithm:ShapeGD}, after a domain deformation, we apply transfer learning 
with training interrupted after at most $200$ steps.

%
%
The Monte Carlo simulations for the components of the probabilistic formulation are accom\-plished with the following samplesizes: we draw $| \samples{\Omega_{\nn}^\pm}{uniform} | = 10^7$ samples to approximate $\lambda^2(\Omega^\pm)$ and the constants \eqref{eqn:acc-ref-constantRandomStart} for the random-start acceptance-rejection sampling. Moreover we use $10^6$ of the samples $\samples{\Omega_{\nn}^\pm}{uniform}$ to approximate the expectation in $\constant{\pm}$; the safety factor is $\safety=1.1$. We draw merely $10^2$ samples of random-starts $\samples{0}{$\pm$}$ and thus obtain the same number of exit samples $\exitSamples$. Due to the low sample size, the resulting Monte Carlo estimate exhibits relatively high variance. However, for our benchmark example the accuracy is sufficient, providing a favorable trade-off between accuracy and computational time.

%
%
We solve the deformation equation using the \textsc{FEniCS} library, see \cite{Fenics_book}, which natively supports finite element function spaces and differential operators on surfaces (our polygonal chain), see \cite{gmd-6-2099-2013}. For the results shown in Figure~\ref{fig:SimulationResult_ATL_ENUMATH}, we operate on $160$ vertices within this polygonal chain. We also use the \textsc{FEniCS} library to visualize the PINN-solution in this figure by computing linear interpolations after evaluating the PINN on a point cloud.

%
%
The objective is evaluated with Monte Carlo using the previously generated $\samples{\Omega_{\nn}^\pm}{uniform}$, i.e.,
\begin{equation}\label{eqn:ObjectiveWithMonteCarlo}
	\begin{aligned}
		\int_{\Omega} \tfrac{1}{2} (u - \trackingData)^2 \de \lambda^2 
		& \approx \lambda^2 (\Omega_{\nn} )\ \E \Big[ \tfrac{1}{2} \big( \nn (\mathcal{U}) - \trackingData (\mathcal{U}) \big)^2 \Big] \\
		& \approx \big( \bar{\prob}_{\text{ref}} [{\Omega^+_{\nn}}] + \bar{\prob}_{\text{ref}} [{\Omega^-_{\nn}}] \big) \sum_{ x\in \samples{\Omega_{\nn}^+}{uniform} \cup \samples{\Omega_{\nn}^-}{uniform}} \tfrac{1}{2} \big( \nn(x) - \trackingData(x)\big)^2,
	\end{aligned}
\end{equation}
where $\mathcal{U}$ denotes a uniformly distributed random variable on $\Omega$, and $\bar{\prob}_{\text{ref}} [{\Omega^\pm_{\nn}}]$ denotes the numerical evaluation of the measures of the subdomains, i.e., the fraction in \eqref{eqn:ApproxMeasurePartition}. 

\section{Perspectives for a Fully Mesh-Free Probabilistic Framework}
The aim of this work is to demonstrate the applicability of a stochastic formulation of the shape derivative. Special attention is placed on the fact, that the probabilistic shape derivative is independent of a given discretization of the shape and can in particular be used in mesh-free approaches. The numerical implementation is heavily based on Monte-Carlo sampling and PINNs, such that the proposed method could potentially circumvent the curse of dimensionality for problems with many spatial dimensions. We use the probabilistic shape derivative to solve a Laplace-type tracking problem, for which we use merely a polygonal chain to keep track of the shape. This is done for convenience and ease of computing the $H^1(\partial \Omega_\text{FEM})$ descent direction, for which mesh-free alternatives also exist.

%
%

We emphasize that the probabilistic shape derivative can be used with any representation of the domain iterate, provided an outer normal $n$ can be defined and uniform sampling from the boundary is possible. Popular choices could be a mesh or a level-set formulation. 

In our application, we train the PINN to solve the PDE on the hold-all $\HoldAll$ with (constant) Dirichlet on an interior interface $\partial \Omega$. Particularly the trained PINN serves as a level-set formulation for the domain, and we can compute $n \approx \frac{\nabla \nn}{\|{\nabla \nn}\|}$.
Furthermore, with the PINN level-set formulation one could generate the required samples to perform the Monte Carlo simulation of the boundary integral in \eqref{eqn:TrackingProblem}, see, e.g., \cite{ling2025uniform, xu2024monte, zappa2018monte}.

%
%
A probabilistic approach to solve the deformation equation could be to compute the $L^2 (\partial\Omega)$ Riesz representant via non-parametric density estimation. For this, let $\nu^\pm$ denote the distribution of the random-start exit-kill variables. Assuming that there are densities $\ell^\pm \colon \partial\Omega \rightarrow [0,\infty)$ such that $\de \nu^\pm = \ell^\pm\de \mathcal{S}^{1}$ the map $\nabla J \defined \trackingData^2 n + \nabla u \big(\constant{+}\ell^+ -\constant{-} \ell^-\big)$ satisfies
\begin{equation}
	\inner{\nabla J}{\dir}_{L^2(\partial \Omega)} \stackrel{!}{=} \D J [\dir] \quad \text{for all admissible directions } \dir
\end{equation}
and thus $\nabla J$ is indeed the $L^2$ shape gradient. A na\"ive estimate could be done using a histogram based on some type of agglomeration on the boundary. For an improved estimation, kernel density estimators on manifolds may be used, see, e.g., \cite{Ozakin2009, Pelletier2005}.

Finally, for transparency we mention our shape optimization loop runs with a fixed stepsize, see Figure~\ref{fig:SimulationResult_ATL_ENUMATH}. Note linesearch strategies typically test for the Armijo--Goldstein condition, i.e., they accept if a step fulfills
\begin{equation*}
	f(x+\alpha p) \leq f(x) + \alpha c \inner{\nabla f(x)}{p}.
\end{equation*}
Here, $\alpha$ is the step length, $p$ is the search direction and $0<c<1$ is a dampening factor. Acceptance of the step is therefore based on the comparison of the actual and expected decrease of the objective. However, in shape optimization, the mesh quality is another limiting factor. Even though our mesh is a mere polygonal line needed for computing the $H^1$-descent, too long a step will still lead to inverted edges, a situation invisible to the Armijo--Goldstein rule. For this reason, we use a manually determined fixed step length well below what would be acceptable by the stricter Armijo--Goldstein condition. Otherwise, one has to base the Armijo criterion on a merit function, which includes mesh quality, see, e.g., \cite[Section~5]{doi:10.1137/24M1646121}.

\bibliographystyle{plain}
\bibliography{references}
\end{document}

%% file: optimization_results_figure.tex
%
%
\begin{figure}[h!]
	\centering
	
	\begin{subfigure}[t]{0.44\linewidth}
		\centering
		\includegraphics[height=6.1cm]{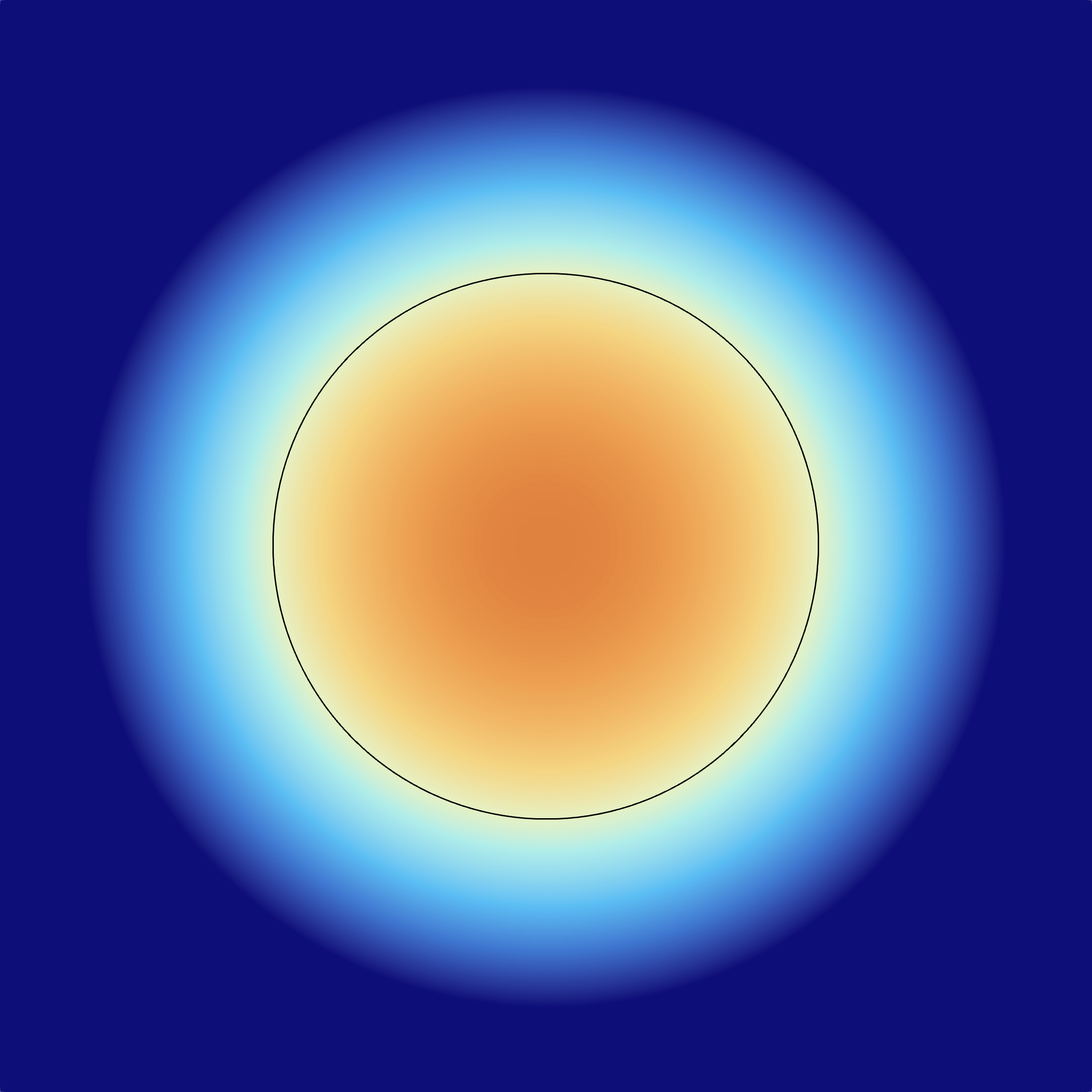}
		\caption{\footnotesize Initial shape. Point evaluation of the PINN with linear interpolation for visualization purposes. Zero level-set $\{ \nn = 0\}$ in black.}
	\end{subfigure}
	\hfill
	\begin{subfigure}[t]{0.54\linewidth}
		\centering
		\includegraphics[height=6.1cm]{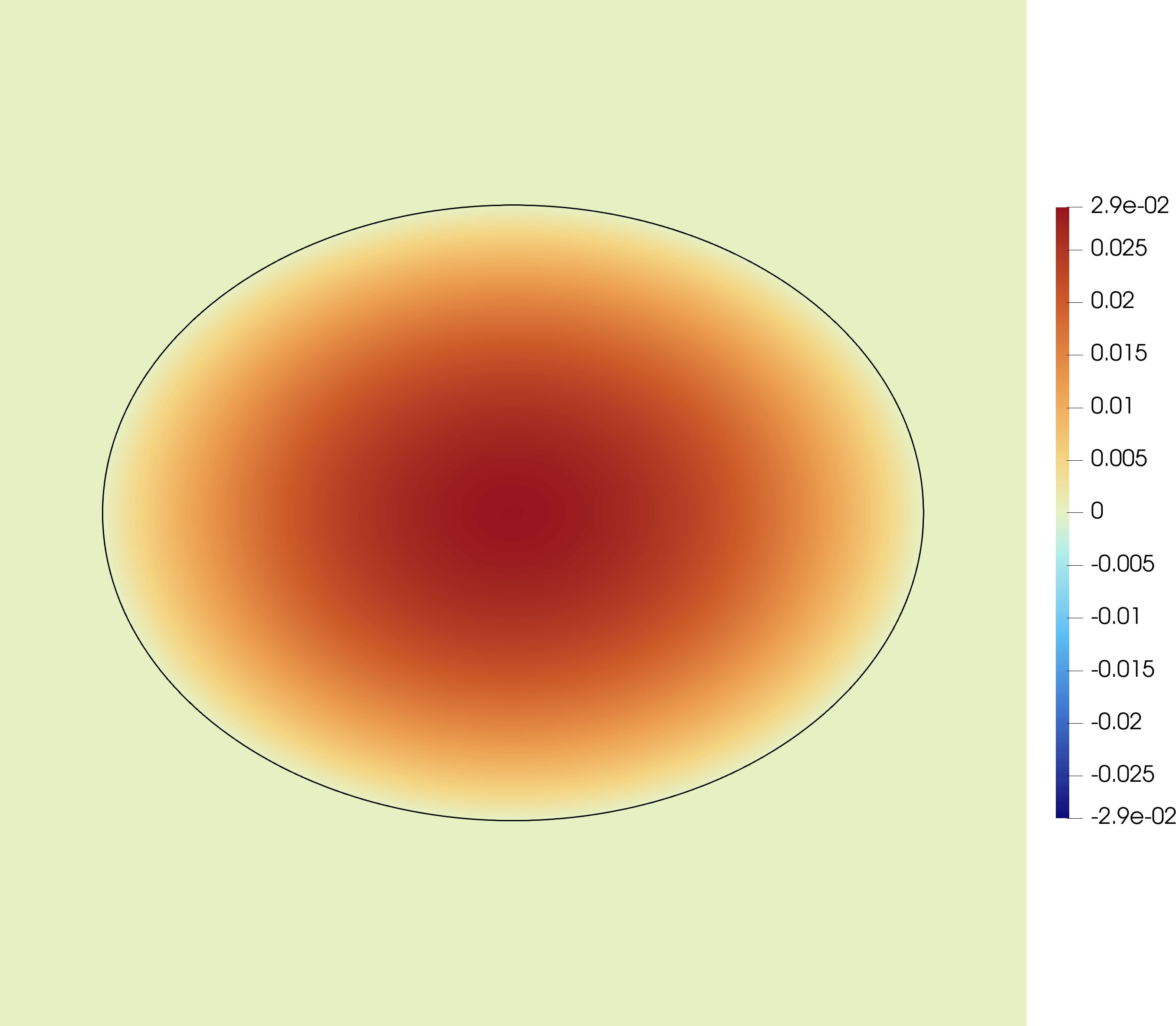}
		\caption{\footnotesize Tracking data $z_d$ on the hold-all $\mathcal{H}$. The Dirichlet conditions holds on the inner interface $\partial \Omega$ and the interior is the analytic solution of the shape optimization problem on $\Omega^*$.}
	\end{subfigure}
	
	\vspace*{0.25cm}
	\begin{subfigure}[t]{0.44\linewidth}
		\centering
		\includegraphics[height=6.1cm]{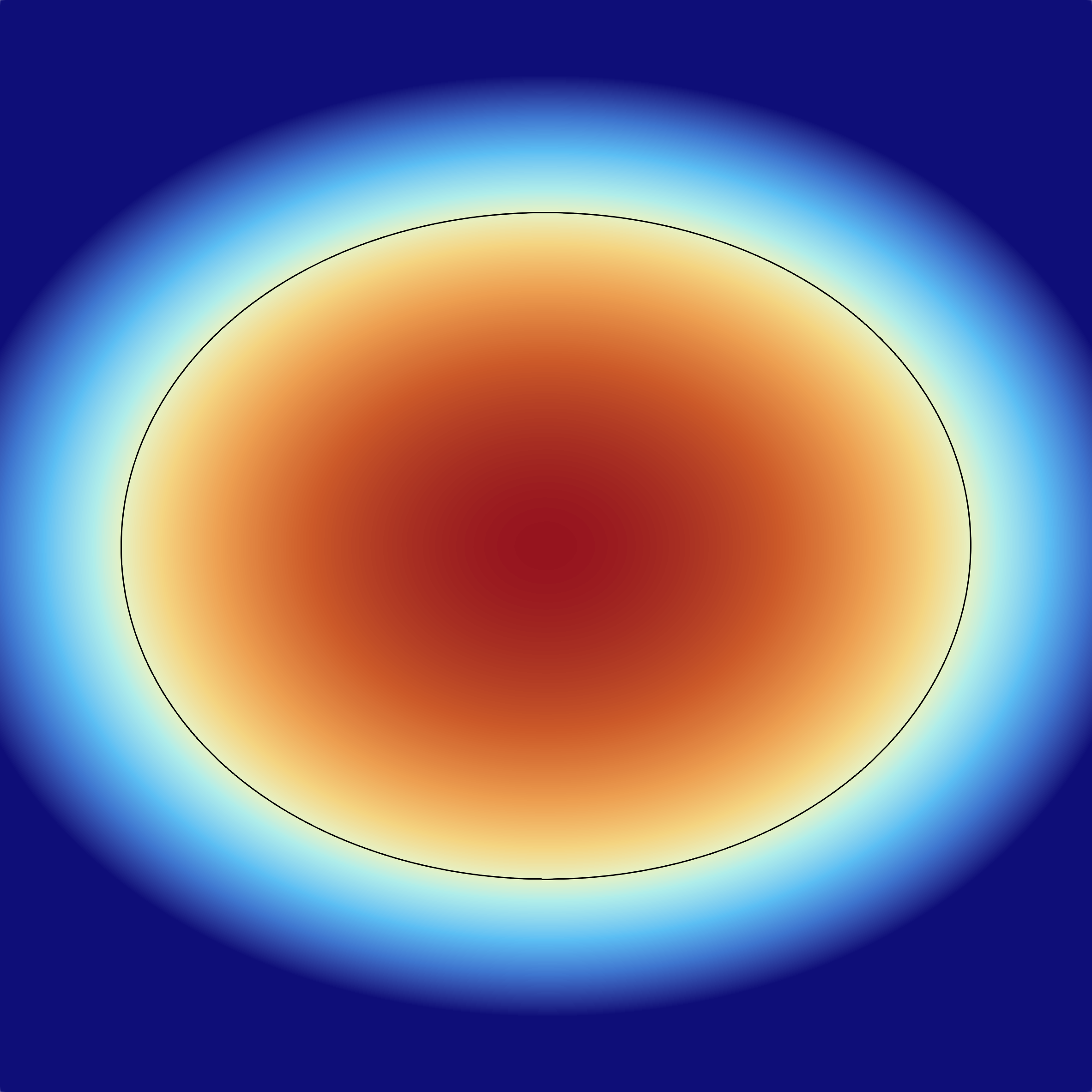}
		\caption{\footnotesize Reconstructed shape $\Omega_{\nn}^*$. Evaluation of the PINN; linearly interpolated for visualization purpose. Level-set $\{ \nn = 0\}$ in black.}
	\end{subfigure}	
	\hfill
	\begin{subfigure}[t]{0.54\linewidth}
		\centering
		\begin{tikzpicture}
			\useasboundingbox (0,0) rectangle (6.1cm, 6.1cm);		
			\begin{axis}[
				at={(0,0.0)}, anchor=south west, width=6.1cm, height=6.1cm,
				xlabel={\footnotesize{Iteration} }, ylabel={\footnotesize{\NoHyper Monte Carlo Objective \eqref{eqn:ObjectiveWithMonteCarlo}\endNoHyper}}, ymode=log, grid=both, scale only axis, 
				y label style={at={(axis description cs: 1.01,0.5)}, anchor=north} 
				]
				\addplot[thick]
				table[x index=0, y index=1] {data/OptHistory_ATL.txt};
			\end{axis}
		\end{tikzpicture}
	\end{subfigure}
	\caption[Optimization Results]{\footnotesize Optimization results.}
	\label{fig:SimulationResult_ATL_ENUMATH}
\end{figure}